\documentclass[12pt]{amsart}
\usepackage{amssymb,amsmath, amsthm,latexsym}
\usepackage{a4wide}

\theoremstyle{plain}

\newtheorem{lem}{Lemma}[section]
\newtheorem{theo}[lem]{Theorem}

\parskip=\medskipamount
\font\k=cmr7

  \newcommand {\loc}{\mbox{\k loc}}

  \newcommand {\C}{{\mathbb C}}
  \newcommand {\Hb}{{\mathbb H}}
  
  \newcommand {\R}{{\mathbb R}}
  \newcommand {\Z}{{\mathbb Z}}
  \newcommand {\Q}{{\mathbb Q}}
  
  \newcommand {\Pf}{{\mathbb P}}

\renewcommand {\H}{{\mathcal H}}

  \newcommand {\cO}{{\mathcal O}}

  \newcommand {\E}{{\mathcal E}}
  \newcommand {\cD}{{\mathcal D}}
 \newcommand {\cP}{{\mathcal P}}
 \newcommand {\cF}{{\mathcal F}}

\newcommand {\ba}{\backslash}

\renewcommand{\Im}{\operatorname{Im}}
\renewcommand{\Re}{\operatorname{Re}}

\newcommand{\Id}{\operatorname{Id}}

\newcommand{\dom}{\operatorname{dom}}

\newcommand{\PSL}{\operatorname{PSL}}

\begin{document}
\title[Spectral interpretation of zeros]
{\Large\bf A spectral interpretation of the\\
 zeros of the constant term of\\ 
\vskip5pt
certain  Eisenstein series}
\date{\today}

\author{Werner M\"uller}
\address{Universit\"at Bonn\\
Mathematisches Institut\\
Beringstrasse 1\\
D -- 53115 Bonn, Germany}
\email{mueller@math.uni-bonn.de}


\begin{abstract}
In this paper we consider the constant term $\varphi_K(y,s)$ of the
non-normalized Eisenstein series attached to $\PSL(2,\cO_K)$, where  $K$ is
either $\Q$ or an imaginary quadratic  field of class number one.
The main purpose of this paper is to show that for every $a\ge 1$ the zeros of
the Dirichlet series $\varphi_K(a,s)$ admit a spectral interpretation  in 
terms of  eigenvalues of a natural 
self-adjoint operator $\Delta_a$. This
implies that, except for at most two real zeros, all zeros of
$\varphi_K(a,s)$ are on the critical
line, and  all zeros are simple. For $K=\Q$ this is due to
Lagarias and Suzuki \cite{LS} and Ki \cite{Ki}.

\end{abstract}
\maketitle

\setcounter{section}{-1}
\section{Introduction}
Let $\Gamma=\PSL(2,\Z)$ be the modular group and $\Hb$ the upper half-plane.
We consider the (completed) 
nonholomorphic Eisenstein series $E^\ast(z,s)$ for the modular
group which is given for $z\in\Hb$ and $\Re(s)>1$ by
\begin{equation}
E^\ast(z,s):=\pi^{-s}\Gamma(s)\left(\frac{1}{2}\sum_{(m,n)\in\Z^2\setminus (0,0)}
\frac{y^s}{|cz+d|^{2s}}\right).
\end{equation}
It is well known that for fixed $z$, $E^\ast(z,s)$ admits a meromorphic 
continuation to the whole $s$-plane, and satisfies the functional equation
\[E^\ast(z,s)=E^\ast(z,1-s).\]
Its only singularities are simple poles at $s=1$ and $s=0$. As a function of
$z$, it is invariant under $\Gamma$
\[E^\ast(\gamma z,s)=E^\ast(z,s), \quad \gamma\in\Gamma.\]
In particular $E^\ast(z,s)$ is invariant under $z\mapsto z+1$ and so it has a
Fourier expansion
\[E^\ast(x+iy,s)=\sum_{n\in\Z}a_n(y,s)e^{2\pi ix}.\]
The zeroth Fourier coefficient $\varphi_0(y,s)$ is given by
\[\varphi_0(y,s)=\Lambda(2s)y^s+\Lambda(2s-1)y^{1-s},\]
where $\Lambda(s)$ is the completed zeta function
\[\Lambda(s)=\pi^{-s/2}\Gamma\left(\frac{s}{2}\right)\zeta(s).\]
Let $a>0$. The zeros of  $\varphi_0(a,s)$ 
have been studied by various authors. Hejhal
\cite[Proposition 5.3(f)]{He} proved that for all $a\ge 1$, the complex zeros
of $\varphi_0(a,s)$ are on the critical line $\Re(s)=1/2$. Lagarias and Suzuki 
\cite{LS} reproved this result and also determined the occurrence of real
zeros. Ki \cite{Ki} proved that all complex zeros are simple.
Putting these results together, we have following theorem.
\begin{theo}\label{th0.1}
For each $a\ge 1$ all complex zeros of $\varphi_0(y,s)$ are simple and lie on
the 
critical line $\Re(s)=1/2$. Moreover there is a critical value
$a^\ast=4\pi e^{-\gamma}=7.055...$
such that the following holds:
\begin{itemize}
\item[1)] For $1\le a\le a^\ast$ all zeros are on the critical line.
\item[2)] For $a>a^\ast$ there are exactly two zeros off the critical
  line. These 
are real simple zeros $\rho_a, 1-\rho_a$ with $\rho_a\in(1/2,1)$. The zero
$\rho_a$ is a nondecreasing function of $a$ and $\lim_{a\to\infty}\rho_a=1$.
\end{itemize}
\end{theo}
The first aim of this paper is to point out that there is actually a
spectral interpretation of the zeros of $\varphi_0(a,s)$, which  
gives another proof of this theorem. For $a> 0$ let $\Delta_a$ be the
cut-off Laplacian $\Delta_a$ introduced by Lax and Phillips \cite{LP}.
It acts in the subspace $\H_a\subset L^2(\Gamma\ba \Hb)$  of all
$f\in L^2(\Gamma\ba \Hb)$ satisfying
$\int_0^1f(x+iy)\,dx=0 $ for almost all $y\ge a$.
The cut-off Laplacian $\Delta_a$ is a nonnegative self-adjoint operator with
pure point spectrum. The spectrum has been studied by Colin de Verdiere 
\cite{CV}. Let 
\begin{equation}\label{0.2}
c(s)=\frac{\Lambda(2s-1)}{\Lambda(2s)},\quad s\in\C.
\end{equation}
The following theorem is an immediate consequence of \cite[Th\'eor\`eme 5]{CV}.

\begin{theo}\label{th0.2}
For every $a>0$, the spectrum of $\Delta_a$ is the union of the cuspidal
eigenvalues 
$0<\lambda_1\le\lambda_2\le\cdots\to\infty$ of $\Delta$ and a sequence of
eigenvalues
$$0<\mu_0(a)< \mu_1(a)<\cdots$$
with the following properties:
\begin{enumerate}
\item[1)] Each eigenvalue $\mu_j(a)$ is a decreasing function of $a$.
\item[2)] If $a\ge1$, each eigenvalue $\mu_j(a)$ has multiplicity 1. Moreover
 $\lim_{a\to\infty}\mu_0(a)=0$ and if $j\ge 1$, then $\mu_j(a)\ge 1/4$ and
 $\lim_{a\to\infty}\lambda_j(a)= 1/4$. 
\item[3)] Let $a\ge1$. Then the  map $s\mapsto s(1-s)$ is a bijection between
  the zeros $\rho\not=1/2$ of $\varphi_0(a,s)$ and the eigenvalues
  $\mu_j(a)\not=1/4$ of $\Delta_a$.
\item[4)] $1/2$ is a zero of $\varphi_0(a,s)$ for all $a>0$.  
There is $j$ with $\mu_j(a)=1/4$ if and only if $c^\prime(1/2)=-2\log a$.
\end{enumerate}
\end{theo}
This theorem implies that for $a\ge 1$  there is at most one
 zero of $\varphi_0(a,s)$ which is off the
line $\Re(s)=1/2$. The simplicity of 
the zeros follows from the consideration of the corresponding eigenfunctions
of $\Delta_a$. 

The main purpose of this article is to extend these results to the constant
term of other Eisenstein series. In the  
present paper we will consider the Eisenstein series attached to
$\PSL(2,\cO_K)$, 
where $\cO_K$ is the ring of integers of
an imaginary quadratic field $K=\Q(\sqrt{-D})$ of class number one, $D$ being
a square free positive integer. These are exactly the fields
$\Q(\sqrt{-D})$ with $D=1,2,3,7,11,19,43,67,163$ \cite{St}.  Let $d_K$ be the
discriminant of $K$ and let $\zeta_K(s)$ be the Dedekind zeta function of
$K$. Let
\begin{equation}
\Lambda_K(s)=\left(\frac{\sqrt{|d_K|}}{2\pi}\right)^s\Gamma(s)\zeta_K(s)
\end{equation}
be the completed zeta function. 
Then $\Lambda_K(s)$ satisfies the functional equation 
$\Lambda_K(s)=\Lambda_K(1-s)$.
For $a>0$ let 
\begin{equation}
\varphi_{K}(a,s):=a^s\Lambda_K(s)+a^{2-s}\Lambda_K(s-1).
\end{equation}
Note that $\varphi_{K}(a,s)$ is a Dirichlet series which satisfies the
functional  equation $\varphi_{K}(a,s)=\varphi_{K}(a,2-s)$.
Let $\xi_K(s)=s(s-1)\Lambda_K(s)$. Then $\xi_K(s)$ is entire. Put
\[a^\ast_K:=\exp\left(1+\frac{\xi_K^\prime(1)}{\xi_K(1)}\right). \]
Then our main result is the following theorem.
\begin{theo}\label{th0.3}
Let $K=\Q(\sqrt{-D})$ be an imaginary quadratic field of class number one. 
Then for each $a\ge 1$ all complex zeros of 
$\varphi_K(a,s)$ are simple and lie on the line $\Re(s)=1$. Moreover 
\begin{enumerate}
\item[1)] For $a>\max\{a^\ast_K,1\}$ there are exactly two zeros off the 
critical line.
These are simple zeros $\rho_a,2-\rho_a$ with $\rho_a\in(1,2)$. The zero
$\rho_a$ is a nondecreasing function of $a$, and $\lim_{a\to\infty}\rho_a=2$. 
\item[2)] If  $1\le a^\ast_K$ and $1\le a\le a^\ast_K$,
 all zeros of $\varphi_K(a,s)$ are on the critical line.
\end{enumerate}
\end{theo}  
To prove Theorem \ref{th0.3}, we extend Colin de Verdiere's Theorem 
\cite[Th\'eor\`eme 5]{CV} to our setting. Let $\cO_K$ be the ring of integers
of $K$ and let $\Gamma=\PSL(2,\cO_K)$. Then $\Gamma$ is a discrete subgroup of 
$\PSL(2,\C)$ which acts properly discontinuously on the 3-dimensional
hyperbolic space $\Hb^3$. The quotient $\Gamma\ba\Hb^3$ has finite volume
and a single cups at $\infty$. The function $\varphi_K(y,s)$ appears to be the 
zeroth Fourier coefficient of the modified Eisenstein series attached to the
cusp $\kappa=\infty$. Let $\Delta_a$ be the corresponding cut-off Laplacian
on $\Gamma\ba\Hb^3$. Then we generalize Theorem \ref{th0.2} to this setting.
As above, this leads to a spectral interpretation of the zeros of
$\varphi_K(a,s)$, and we deduce Theorem \ref{th0.3} from this spectral
interpretation. 

The constant $a^\ast_K$ can be computed using the Kronecker limit formula 
\cite[p. 273]{La} 
and the Chowla-Selberg formula \cite[(2), p.110]{SC}.
For example, we get
\[a^\ast_{\Q(i)}=\frac{4\pi^2 e^{2+\gamma}}{\Gamma(1/4)^4}\approx 3.00681,
\quad
a^\ast_{\Q(\sqrt{-2})}=\frac{8\pi^2e^{2+\gamma}}
{\left(\Gamma\left(\frac{1}{8}\right)\Gamma\left(\frac{3}{8}\right)\right)^{2}}
\approx 3.2581 .\] 
Thus $a^\ast_{\Q(\sqrt{-D})}>1$ for $D=1,2$ and therefore, in the range 
$1\le a\le a^\ast_{\Q(\sqrt{-D})}$  all zeros of
$\varphi_{\Q(\sqrt{-D})}(a,s)$ are  
on the line $\Re(s)=1$. 

The method used to prove Theorem \ref{th0.3} can be extended so that
an arbitrary number field $K$ of class number one can be treated. 
The underlying global Riemannian symmetric 
space is $X=\Hb^{r_1}\times (\Hb^3)^{r_2}$, where $r_1$ is the number of real
places and $r_2$ the number of pairs of complex conjugate places of $K$. The
structure of $\Gamma\ba X$ is slightly more complicated, however
the proof of the corresponding statement is completely analogous. More
generally, it seems to be possible to prove similar results for the
constant terms of rank one cuspidal Eisenstein series attached to
$\PSL(n,\Z)$.  
 
\section{The cut-off Laplacian}
\setcounter{equation}{0}

The cut-off Laplacian can be defined for every rank one locally symmetric 
space. In the present paper we are dealing only with hyperbolic
manifolds of dimension 2 and 3. So  we discuss only the case of a
hyperbolic 3-manifold. The general case is similar.

Let  
\[\Hb^3=\{(x_1,x_{2},y)\in\R^3\colon y>0\}\]
be the hyperbolic 3-space with its hyperbolic metric
\[ds^2=\frac{dx_1^2+dx_{2}^2+dy^2}{y^2}.\]
The hyperbolic Laplacian $\Delta$ is given by
\[\Delta=-y^2\left(\frac{\partial^2}{\partial x_1^2}
+\frac{\partial^2}{\partial x_2^2}+ 
\frac{\partial^2}{\partial y^2}\right)+y\frac{\partial}{\partial y}.\]
If we regard $\Hb^3$ as the set of all quaternions $z+yj$ with $z\in\C$
and $y>0$, then $G=\PSL(2,\C)$ is the group of all orientation preserving
isometries of $\Hb^3$. It acts by linear fractional transformations, i.e., for
$\gamma=\begin{pmatrix}a&b\\c&d\end{pmatrix}\in G$,
\[\gamma(w)=(aw+b)(cw+d)^{-1},\quad w\in\Hb^3.\]
Let
$$B(\C):=\left\{\begin{pmatrix}\lambda&z\\0&\lambda^{-1}\end{pmatrix}\colon
  \lambda\in\C^\ast,\; z\in\C\right\}/\{\pm \Id\},\quad
N(\C):=\left\{\begin{pmatrix}1&z\\0&1\end{pmatrix}\colon z\in\C\right\}.$$ 
Let $\Gamma\subset G$ be a discrete subgroup of finite co-volume. 
Let $\kappa_1=\infty, \kappa_2,...,\kappa_m\in \Pf^1(\C)$ be a complete set of
$\Gamma$-inequivalent cusps, and let $\Gamma_i$ be the stabilizer of
$\kappa_i$ in $\Gamma$. Choose $\sigma_i\in G$ such that
$\sigma_i(\kappa_i)=\infty$, $i=1,...,m$. Then
\begin{equation}\label{1.1}
\sigma_i\Gamma_i\sigma_i^{-1}\cap
N(\C)=\left\{\begin{pmatrix}1&\lambda\\0&1\end{pmatrix}\colon 
\lambda\in L_i\right\},
\end{equation}
where $L_i$ is a lattice in $\C$ \cite[Theorem 2.1.8]{EGM}. Note that
for $D\not=1,3$  the intersection actually coincides with  
$\sigma_i\Gamma_i\sigma_i^{-1}$. 
Choose closed fundamental sets $\cP_i$  for
the action of $\sigma_i\Gamma_i \sigma_i^{-1}$ 
on $\Pf^1(\C)\setminus\{\infty\}=\C$. For $T>0$ define
\[\tilde \cF_i(T):=\left\{ (z,y)\in\Hb^3\colon z\in\cP_i,\;\; y\ge T\right\}.\] 
Let $\cF_i(T)=\sigma_i^{-1}(\tilde \cF_i(T))$. There exists $T_i>0$ such 
that
any two points in $\cF_i(T_i)$ are $\Gamma$-equivalent if and only if they are
$\Gamma_i$-equivalent. For such a choice of $T_1,...,T_m$  there exists a
compact subset $\cF_0\subset \Hb^3$ such that
\begin{equation}\label{1.2}
\cF:=\cF_0\cup \cF_1(T_1)\cup\cdots\cup \cF_m(T_m)
\end{equation}
is a fundamental domain for $\Gamma$ \cite[Proposition 2.3.9]{EGM}. Let
\begin{equation}\label{1.3}
b:=\max\{T_1,...,T_m\}.
\end{equation}
By $C^\infty_c(\Gamma\ba\Hb^3)$ we denote
the space of $\Gamma$-invariant $C^\infty$-functions on $\Hb^3$ with compact 
support in $\cF$. For $f\in C^\infty_c(\Gamma\ba \Hb^3)$ let
$$\parallel f\parallel^2=\int_{\cF}|f(x_1,x_2,y)|^2\;\frac{dx_1dx_2dy}{y^3},$$
and let $L^2(\Gamma\ba \Hb^3)$ be the completion of $C^\infty_c(\Gamma\ba
\Hb^3)$ with respect to this norm.
Similarly let
$$\parallel df\parallel^2=\int_{\cF} df\wedge\ast \overline{df}
=\int_{\cF}\parallel
df(x_1,x_2,y)\parallel^2\;\frac{dx_1dx_2dy}{y}.$$ 
Let $H^1(\Gamma\ba \Hb^3)$ denote the completion of $C^\infty_c(\Gamma\ba
\Hb^3)$ with respect to the norm
$$\parallel f\parallel^2_1:=\parallel f\parallel^2+\parallel df\parallel^2.$$
Denote by $|\cP_i|$ the Euclidean area of the fundamental domain 
$\cP_i\subset\C$ of
$\sigma_i\Gamma_i\sigma_i^{-1}$.
For $f\in L^2(\Gamma\ba \Hb^3)$ let 
\[f_{j,0}(y)=\frac{1}{|\cP_i|}\int_{\cP_i} f(\sigma_i^{-1}(x_1,x_2,y))\;dx_1dx_2\]
be the zeroth coefficient of the  Fourier expansion of $f$ at the cusp
$\kappa_j$. Note that for $f\in H^1(\Gamma\ba \Hb^3)$, each $f_{j,0}$ belongs to
$H^1(\R^+)$ and therefore, each $f_{j,0}$ is a continuous function on $\R^+$. 
For $a>0$ let
\[H^1_a(\Gamma\ba \Hb^3):=\left\{f\in H^1(\Gamma\ba \Hb^3)\colon f_{j,0}(y)=0\;
\;{\mathrm for}\;y\ge a,\;j=1,...,m\right\}.\]
Then $H^1_a(\Gamma\ba\Hb^3)$ is a closed subspace of $H^1(\Gamma\ba\Hb^3)$. 
Hence the quadratic form 
\[q_a(f)=\parallel df\parallel^2,\quad f\in H^1_a(\Gamma\ba\Hb^3),\]
is closed. Let $\Delta_a$ denote the self-adjoint operator which represents
the quadratic form $q_a$. It acts in the Hilbert space $\H_a$ which is the
closure of $H^1_a(\Gamma\ba\Hb^3)$ in $L^2(\Gamma\ba\Hb^3)$. By definition,
$\Delta_a$ is nonnegative. Its domain can be described as follows.
Let $\psi_{j,a}\in\cD^\prime(\Gamma\ba\Hb^3)$ be defined by 
\[\psi_{j,a}(f):=f_{j,0}(a),\quad f\in C^\infty_c(\Gamma\ba\Hb^3).\]
Let $b>0$ be defined by (\ref{1.3}). 
\begin{lem} Let $a\ge b$. Then the domain of $\Delta_a$ consists of all 
$f\in H^1_a(\Gamma\ba \Hb^3)$ for which there exist $C_1,...,C_m\in\C$ 
such that
\begin{equation}\label{1.4}
\Delta f-\sum_{j=1}^mC_j\psi_{j,a}\in L^2(\Gamma\ba \Hb^3).
\end{equation}
\end{lem}
The proof of the lemma is analogous to the proof of Th\'eor\`em 1 in
\cite{CV}. 
Let  $f\in \dom(\Delta_a)$. By the lemma, there exist $C_1,...,C_m\in\C$ such 
that (\ref{1.4}) holds. Then $\Delta_af$ is given by
\begin{equation}\label{1.5}
\Delta_af=\Delta f-\sum_{j=1}^mC_j\psi_{j,a}.
\end{equation}

Furthermore, we have
\begin{lem}\label{l1.2}
$\Delta_a$ has a compact resolvent.
\end{lem}
\begin{proof} The proof is a simple extension of the argument in 
\cite[p.206]{LP}. 
\end{proof}
So the spectrum of $\Delta_a$  consists of a sequence of eigenvalues 
$0\le \lambda_1(a)\le \lambda_2(a)\le\cdots\to\infty$ with finite
multiplicities. To describe the eigenvalues and eigenfunctions of $\Delta_a$
more explicitely, we need  to consider the Eisenstein series. Let $\sigma_j\in
G$, $j=1,...,m$, be as above. 
For each cusp $\kappa_i$,  the Eisenstein series $E_i(w,s)$ attached to
$\kappa_i$ is defined to be
\[E_i(w,s):=\sum_{\gamma\in\Gamma_i\ba \Gamma}y(\sigma_i\gamma w)^s,\quad
\Re(s)>2,\]
where $y(\sigma_i\gamma w)$ denotes the $y$-component of $\sigma_i\gamma w$. 
Selberg has shown that $E_i(w,s)$ can be meromorphically continued in $s$ to 
the whole complex plane $\C$, and it is an automorphic eigenfunction of
$\Delta$ with 
$$\Delta E_i(w,s)=s(2-s)E_i(w,s).$$
Since $E_i(w,s)$ is $\Gamma$-invariant, it is invariant under the stabilizer
$\Gamma_j$ of the cusp $\kappa_j$ and therefore, admits a Fourier expansion
at  $\kappa_j$. The constant Fourier coefficient is given by
\begin{equation}\label{1.3a}
\delta_{ij}y^s+C_{ij}(s)y^{2-s}.
\end{equation}
Put
\[C(s):=\left(C_{ij}(s)\right)_{i,j=1,...,m}.\]
This is the so called ``scattering matrix''. The scattering matrix and the 
Eisenstein series satisfy the following system of functional equations.
\begin{equation}\label{1.7}
\begin{split}
&C(s)C(2-s)=\Id,\\
&E_i(w,s)=\sum_{j=1}^m C_{ij}(s)E_j(w,2-s), \quad i=1,...,m.
\end{split}
\end{equation}

Now recall that a square integrable eigenfunction $f$ of $\Delta$ is called 
cuspidal, if the zeroth Fourier coefficient $f_{j,0}$ of $f$
at the cusp $\kappa_j$ vanishes for all $j=1,...,m$. Denote by
$S(\lambda;\Gamma)$ the space of cuspidal eigenfunctions of $\Delta$ with 
eigenvalue $\lambda$. A function $f\in C^0(\Gamma\ba \Hb^3)$ is called of 
moderate growth, if there exists $R\in\R^+$ such that
\[f(\sigma^{-1}_i(x_1,x_2,y))=O(y^R),\quad {\mathrm
  for}\;\;(x_1,x_2,y)\in\cF_i(b),\; i=1,...,m.\] 
It follows from the Fourier expansion in the cusps \cite[Section 3.3.3]{EGM}
that a  cuspidal eigenfunction $f$ of $\Delta$ is rapidly  decreasing in each
cusp. Therefore for $\psi\in\E(\lambda,\Gamma)$ and $\varphi\in
S(\lambda,\Gamma)$ , the inner product $\langle \psi, \varphi\rangle$ 
is well defined. For $\lambda\in\C$ let
\begin{equation}
\E(\lambda;\Gamma):=\left\{\psi\in L^2_{\loc}(\Gamma\ba\Hb^3)\colon
  \Delta\psi=\lambda\psi,\;\psi\;\;{\mathrm is \; of\;
  moderate\; growth},\;\psi\perp S(\lambda;\Gamma)\right\}.
\end{equation}
\begin{lem}\label{l1.3} Let $m$ be the number of cusps of $\Gamma\ba\Hb^3$. 
Then 
\begin{enumerate}
\item[1)] For every $\lambda\in\C$ we have $\dim\E(\lambda;\Gamma)\le m$. 
\item[2)] Suppose that $\lambda=s(2-s)$, where $\Re(s)\ge 1$, $s\not=1$,
  and $s$ is not a pole of any Eisenstein series. Then
$E_1(w,s),...,E_m(w,s)$ is a basis of $\E(\lambda;\Gamma)$.
\end{enumerate}
\end{lem}
\begin{proof}
The proof is analogous to the proof of Satz 11 in \cite[p.171]{Ma}. Let 
$\phi,\psi\in\E(\lambda;\Gamma)$ and $\lambda=s(2-s)$. The zeroth Fourier 
coefficient of $\phi$, $\psi$ at the $j$-th cusp is given by
\[\phi_{0,j}(y)=a_jy^s+b_jy^{2-s},\quad
\psi_{0,j}(y)=c_jy^s+d_jy^{2-s},\;s\not=1,\] 
and for $s=1$ by
\[\phi_{0,j}(y)=a_jy+b_jy\log y,\quad \psi_{0,j}(y)=c_jy+d_jy\log y.\]
Let $\chi_{j,[a,\infty)}$ denote the characteristic function of $\Gamma_i\ba
\cF_i(a)$ in $\Gamma\ba\Hb^3$. Let
\[\phi_a=\phi-\sum_{j=1}^m\chi_{j,[a,\infty)}\phi_{0,j},\quad  
\psi_a=\psi- \sum_{j=1}^m\chi_{j,[a,\infty)}\psi_{0,j}\] 
be the truncations of $\phi$ and $\psi$, respectively, at level $a\ge b$. 
Let $s\not=1$. Integrating by parts, we get
\begin{equation}\label{1.9}
\begin{split}
0=\int_{\Gamma\ba \Hb^3}(\phi_a\Delta\psi_a-\psi_a\Delta\phi_a)\;d{\mathrm 
vol}&=
\sum_{j=1}^m\left(\psi_{0,j}(y)\frac{\partial\phi_{0,j}}{\partial y}(y)-
\phi_{0,j}(y)\frac{\partial\psi_{0,j}}{\partial y}(y)\right)\bigg|_{y=a}\\
&=(2-2s)\sum_{j=1}^m(a_jd_j-b_jc_j).
\end{split}
\end{equation}
For $s=1$ the right hand side equals $\sum_{j=1}^m(a_jd_j-b_jc_j)$.

An element of $\E(\lambda)$ is uniquely determined by its zeroth Fourier 
coefficients. Define $i\colon \E(\lambda)\to \C^{2m}$ by
$i(\phi):=(a_1,b_1,...,a_m,b_m),$
where $a_j,b_j$ are the $0$-th Fourier coefficients of $\phi$ at the $j$-th
cusp. Then $i$ is an 
embedding. Let 
\[I(x,y):=\sum_{j=1}^m(x_{2j-1}y_{2j}-x_{2j}y_{2j-1})\]
be the standard symplectic form on $\C^{2m}$. Then by (\ref{1.9}) and the
corresponding statement for $s=1$,
$\E(\lambda)$ is an isotropic subspace for $I$. Hence $\dim\E(\lambda)\le m$. 

If $s\not=1$ and $s$ is not a pole of $E_j(w,s)$, $j=1,...,m$, the Eisenstein 
series have linearly independent $0$-th Fourier coefficients given by 
(\ref{1.3a}). This implies that $E_1(w,s),...,E_m(w,s)$ are linearly
independent and hence, form a basis of $\E(\lambda,\Gamma)$, where
$\lambda=s(2-s)$.
\end{proof}
We are now ready to describe the spectrum of $\Delta_a$. Let
$0=\Lambda_0<\Lambda_1<\cdots<\Lambda_N<1$ be the eigenvalues of $\Delta$ with
non-cuspidal eigenfunctions. Then $\Lambda_j$ has the form 
$\Lambda_j=\sigma_j(2-\sigma_j)$ with $\sigma_j\in (1,2]$ and $\sigma_j$ is a 
simple pole of some 
Eisenstein series $E_k(w,s)$. The corresponding residue is an 
eigenfunction of $\Delta$ with eigenvalue $\Lambda_j$ 
\cite[Proposition 6.2.2]{EGM}.
The following theorem
is an extension of Theorem 5 in \cite{CV} to the 3-dimensional case.
\begin{theo}\label{th1.4}
The spectrum of $\Delta_a$ is the union of the cuspidal 
eigenvalues $(\lambda_i)_i $ of $\Delta$ and a sequence
$0<\mu_0(a)\le \mu_1(a)\le \cdots $
with the following properties:
\begin{enumerate}
\item[1)] Each $\mu_j(a)$ is a decreasing function of $a$,
\item[2)] Let $a\ge b$. Then each $\mu_j(a)$ has multiplicity $\le m$, and
 the map $s\mapsto s(2-s)$ is a bijection between the
eigenvalues $\mu_j(a)\not\in\{1,\Lambda_0,...,\Lambda_N\}$ and the zeros of
$\rho\in\C\setminus\{1,\sigma_0,...,\sigma_N\}$ of 
$\varphi(s):=\det(C(s)+a^{2s-2}\Id)$.
\end{enumerate}
\end{theo}
\begin{proof}
Let $f\in L^2(\Gamma\ba\Hb^3)$ be a cuspidal eigenfunction of $\Delta$.
Then $f$ belongs to
$H^1_a(\Gamma\ba \Hb^3)$ and $\Delta f$ is square integrable.  Hence
by (\ref{1.4}) and (\ref{1.5}), $f\in\dom(\Delta_a)$ and $\Delta_af=\Delta f$.
Thus all cuspidal eigenfunctions of $\Delta$ are also eigenfunctions of
$\Delta_a$ with the same eigenvalues. 

For 1)
observe that for $a\le a^\prime$ we have
\[H^1_a(\Gamma \ba\Hb^3)\subset H^1_{a^\prime}(\Gamma \ba\Hb^3).\]
Since $\Delta_a$ is a non-negative self-adjoint operator with compact
resolvent, 1) follows from the mini-max characterization of its eigenvalues.

2) Let $\Phi$ be an eigenfunction of $\Delta_a$ with eigenvalue
$\mu_j(a)=s_j(2-s_j)$ and suppose that $\Phi\perp
S(\mu_j(a);\Gamma)$.  Let $\Phi_{i,0}$ be the zeroth Fourier coefficient of
$\Phi$ at the cusp $\kappa_i$. There exists $1\le j\le m$ such that
$\Phi_{j,0}\not\equiv 0$. Let $\hat\Phi_{i,0}$ denote the analytic continuation
of $\Phi_{i,0}$ from $(0,a)$ to $\R^+$. Put
\[\hat\Phi(w):=\Phi(w)+\sum_{i=1}^m\chi_{i,[a,\infty)}\hat\Phi_{i,0}.\]
Then $\hat\Phi\in C^\infty(\Gamma\ba \Hb^3)$ and $\Delta \hat\Phi
=\mu_j(a)\hat\Phi$. 
Moreover $\hat\Phi$
is of moderate growth and $\hat\Phi\perp S(\mu_j(a),\Gamma)$. Therefore $\hat\Phi\in\E(\mu_j(a);\Gamma)$. 
By Lemma \ref{l1.3},
1), it follows that the multiplicity of $\mu_j(a)$ is bounded by $m$.

Next suppose that
$\mu_j(a)=s_j(2-s_j)\not\in\{1,\Lambda_0,...,\Lambda_N\}$. We can  
assume that $\Re(s_j)\ge 1$. As explained above, $s_j$ is not a pole of any
Eisenstein series $E_k(w,s)$, $k=1,...,m$. Therefore, by Lemma \ref{l1.3} 
there exist 
$c_1,...,c_m\in\C$ such that 
\[\hat\Phi(w)=\sum_{i=1}^m c_i E_i(w,s_j).\]
By construction, the zeroth Fourier coefficient $\hat\Phi_{l,0}(y)$ of 
$\hat\Phi(w)$ at
$\kappa_l$ vanishes at $y=a$ for all $l=1,...,m$. By (\ref{1.3a}) this implies
\[\sum_{i=1}^mc_i\left(\delta_{il}a^{s_j}+C_{il}(s_j)a^{2-s_j}\right)=0,\quad
l=1,...,m.\]
Hence $\det(C(s_j)+a^{2s_j-2}\Id)=0$.

Now assume that $\det(C(s_j)+a^{2s_j-2}\Id)=0$, $\Re(s_j)\ge 1$,
 and $s_j\not\in\{1,\sigma_0,...,\sigma_N\}$. Then $s_j$ is not a pole of any
Eisenstein series. Let $0\not=\psi\in\C^m$ such
that
\begin{equation}\label{1.10}
C(s_j)\psi=-a^{2s_j-2}\psi.
\end{equation}
Let $\psi=(c_1,...,c_m)$ and set
\[\hat\Phi(w):=\sum_{k=1}^m c_kE_k(w,s_j).\] 
By (\ref{1.3a}) and (\ref{1.10}), the $0$-th Fourier coefficient of 
$\hat\Phi(w)$ at the cusp $\kappa_l$ is given by
\begin{equation}
\sum_{k=1}^m
c_k\left(\delta_{kl}y^{s_j}+C_{kl}(s_j)y^{2-s_j}\right)
=c_l\left(y^{s_j}-a^{2s_j-2}y^{2-s_j}\right). 
\end{equation}
Therefore, the $0$-th Fourier coefficient $\hat\Phi_{l,0}(y)$ vanishes at
$y=a$ for all $l=1,...,m$. Put
\[\Phi(w)=\hat\Phi(w)-\sum_{i=1}^m\chi_{i,[a,\infty)}\hat\Phi_{l,0}(y).\]
Then $\Phi\in H^1_a(\Gamma\ba\Hb)$ and it follows from the description of
the domain of $\Delta_a$ that $\Phi\in\dom(\Delta_a)$ and 
$\Delta_a \Phi=s_j(2-s_j)\Phi$. 
\end{proof}

Now assume that $\Gamma\ba \Hb^3$ has a single cusp $\kappa=\infty$. Then
there is a single  
Eisenstein series $E(w,s)$, which is given by
\[E(w,s)=\sum_{\gamma\in\Gamma_\infty\ba\Gamma}y(\gamma w)^s,\quad \Re(s)>2.\]
The zeroth Fourier coefficient of $E(w,s)$ at $\infty$ equals
\begin{equation}
\varphi_0(y,s):=y^s+c(s)y^{2-s},
\end{equation}
where $c(s)$ is a meromorphic function of $s\in\C$. 
The following facts are well known 
\cite[pp.243-245]{EGM} .
The poles of $E(w,s)$ in the half-plane $\Re(s)>1$ are 
contained in the interval $(1,2]$ and are all simple. The residue of $E(w,s)$
at a pole 
$\sigma\in (1,2]$ is a square integrable eigenfunction of $\Delta$ with 
eigenvalue $\sigma(2-\sigma)$, which  is non-cuspidal. Moreover all
non-cuspidal eigenfunctions of $\Delta$ are obtained in this way. Thus
the non-cuspidal eigenvalues have multiplicity one. So in this case we get 
the following refinement of Theorem \ref{th1.4}.
\begin{theo}\label{th1.5}
Assume that $\Gamma\ba\Hb^3$ has a single cusp.
Then the spectrum of $\Delta_a$ is the union of the cuspidal 
eigenvalues $(\lambda_i)_i $ of $\Delta$ and a sequence
$0<\mu_0(a)< \mu_1(a)< \cdots $
with the following properties:
\begin{enumerate}
\item[1)] Each $\mu_j(a)$ is a decreasing function of $a$;
\item[2)] Let $a\ge b$. Then each $\mu_j(a)$ has multiplicity 1 and the map
$s\mapsto s(2-s)$ is a bijection between the zeros $\rho\not=1$ of 
$\varphi_0(a,s)$ and the eigenvalues $\mu_j(a)\not=1$ of $\Delta_a$. 
\item[3)] There is $j$ with $\mu_j(a)=1$  if and only if $c(1)=-1$ and
 $c^\prime(1)=-2\log a$.  
\item[4)] Let $0=\Lambda_0<\Lambda_1<\cdots<\Lambda_N<1$ be the
 eigenvalues of $\Delta$ 
with  non-cuspidal eigenfunctions. If $a\ge b$, then 
\[0=\Lambda_0<\mu_0(a)<\Lambda_1<\mu_1(a)<\cdots <\Lambda_N< 
\mu_N(a).\]
Moreover, 
\begin{equation}
\lim_{a\to\infty}\mu_j(a)=\begin{cases}\Lambda_j&,\; 0\le j\le N;\\
1&, \;j>N.
\end{cases}
\end{equation}
\end{enumerate}
\end{theo}
\begin{proof}
1) follows from Theorem \ref{th1.4}. If $s\not\in\{1,\sigma_0,...,\sigma_N\}$,
then 2) follows also from 2) of Theorem \ref{th1.4}. Now suppose that $s_0$
is a pole of $E(w,s)$ in the half-plane $\Re(s)\ge 1$. Then the residue
$\psi$ of $E(w,s)$ at $s_0$ is an eigenfunction of $\Delta$ with non-vanishing
$0$-th Fourier coefficient. Moreover $\psi\perp S(\lambda,\Gamma)$.
Therefore $\psi\in\E(\lambda,\Gamma)$, where 
$\lambda=s_0(2-s_0)$. By Lemma \ref{l1.3} it follows that
 $\dim\E(\lambda,\Gamma)=1$.
Moreover the constant term $\psi_0(y)$ of $\psi$ has the form 
$\psi_0(y)=cy^{2-s_0}$. Especially, it never vanishes on $\R^+$. On the other
hand, if $\Phi$ were an eigenfunction of $\Delta_a$ with eigenvalue $\lambda$,
then the corresponding eigenfunction $\hat\Phi\in\E(\lambda,\Gamma)$, 
constructed in the proof of Theorem \ref{th1.4} would have a constant term
which  
vanishes at $y=a$. Hence the eigenvalues $\Lambda_j$ can not be eigenvalues of
$\Delta_a$. This proves 2)

For 3) we note that $c(1)^2=1$ by the
functional equation (\ref{1.7}). Thus $c(1)=\pm 1$. If $c(1)=-1$,  then 
$E(w,1)\equiv 0$. Put
$\psi(w):=\frac{d}{ds}E(w,s)\big|_{s=1}.$ Then $\Delta \psi=\psi$ and the 
zeroth Fourier coefficient of $\psi$ is given by
\[\psi_0(y)=(2\log y+c^\prime(1))y.\]
Suppose that $c^\prime(1)=-2\log a$. Put 
$\hat \psi_a(w):=\psi(w)-\chi_{[a,\infty)}\psi_0(y(w))$. Then $\hat \psi_a$
is in the domain of $\Delta_a$ and $\Delta_a\hat\psi_a=\hat\psi_a$. 

For the other direction suppose that $\hat\psi$ is an eigenfunction of
$\Delta_a$ with eigenvalue 
$1$ and $\hat\psi\perp S(1,\Gamma)$. Let $\hat\psi_0(y)$ be the $0$-th Fourier
coefficient of $\hat\psi$.  
Extend $\hat\psi_0(y)$ in the obvious way  from $(b,a)$  to a smooth function 
$\psi_0$ on $(b,\infty)$. Set
\[\psi:=\hat\psi+\chi_[a,\infty)\psi_{0}.\]
Then $\psi$ is smooth, of moderate growth, satisfies $\Delta\psi=\psi$ and
$\psi\perp S(1,\Gamma)$. Thus $\psi\in\E(1,\Gamma)$. Therefore by Lemma
\ref{l1.3} it follows that $\dim\E(1,\Gamma)=1$. On the other hand, we have
$0\not=\frac{d}{ds}E(s,w)\big|_{s=1}\in\E(1,\Gamma)$.  
Therefore it follows that there exists
$c\not=0$ such that $\psi(w)=c\frac{d}{ds}E(w,s)\big|_{s=1}$. 
Comparing the constant terms, it follows that $c^\prime(1)=-2\log a$. 

4) follows from the mini-max principle and the fact that $\cup_{a\ge b}
H^1_a(\Gamma \ba\Hb^3)$ is dense in $H^1(\Gamma \ba\Hb^3)$.
\end{proof}

\section{Hyperbolic surfaces of finite area}
\setcounter{equation}{0}

In this section we prove Theorem \ref{th0.2} and deduce Theorem
\ref{th0.1} from it.

Let $\Gamma=\PSL(2,\Z)$ be the modular group. Then $\Gamma\ba \Hb$ has a single
cusp $\kappa=\infty$. As fundamental domain we take the standard domain
\[F:=\left\{z\in\Hb\colon |\Re(z)|<1/2,\,|z|>1\right\}.\]
So we can take $b=1$. 
The Eisenstein series attached to the cusp $\infty$ is given by
\[E(z,s)=\sum_{\gamma\in\Gamma_\infty\ba\Gamma}\Im(\gamma z)^s=\sum_{(m,n)=1}
\frac{y^s}{|mz+n|^{2s}},\quad \Re(s)>1.\]
The constant term of $E(z,s)$ equals
\begin{equation}\label{2.1}
y^s+c(s)y^{1-s},
\end{equation}
where $c(s)$ is the meromorphic function defined by (\ref{0.2}). 

\noindent
{\bf Proof of Theorem 0.2:}
Since the 
the zeros of the completed zeta function $\Lambda(s)$ are all contained in the
strip $0<\Re(s)<1$, it follows that $\Lambda(2s-1)$ and $\Lambda(2s)$ have no
common  zeros. This implies that the zeros of 
\[\varphi_0(a,s)=\Lambda(2s)a^s+\Lambda(2s-1)a^{1-s}\]
coincide with the zeros of
$a^s+c(s)a^{1-s}$. Now note that $\Lambda(s)$ has poles of order 1 at $s=1,0$.
The residue at $s=1$ is 1 and the residue at $s=0$ is $-1$. Hence
\begin{equation}\label{2.2}
c(1/2)=\lim_{s\to 1/2}\frac{\Lambda(2s-1)}{\Lambda(2s)}=-1.
\end{equation}
Now the statements 1), 3) and 4) follow immediately from \cite[Th\'eor\`em
5]{CV}.  For 2) we use that by Roelcke \cite{Ro}, the smallest positive
eigenvalue $\lambda_1$ of $\Delta$ on $\Gamma\ba \Hb$ satisfies
$\lambda_1>1/4$. Then 2) follows from Th\'eor\`em 5, (iii), of \cite{CV}. 

\noindent
{\bf Proof of Theorem 0.1:} Let $a\ge 1$ and let $\rho\not=1/2$ be a zero of
$\varphi_0(a,s)$. Then by Theorem 0.2, 3), $\rho(1-\rho)$ is an eigenvalue of
$\Delta_a$. Hence $\rho(1-\rho)$ is a non-negative real number. If 
$\rho(1-\rho)> 1/4$, then $\rho$ is a complex zero with $\Re(\rho)=1/2$. 
By Theorem 0.2, 2), there exist at most two zeros, $\rho_a$ with 
$0<\rho_a(1-\rho_a)<1/4$. Thus $\rho_a$ and $1-\rho_a$ are zeros and we
may assume that 
$\rho_a\in (1/2,1)$. Let $a<a^\prime$. Then by 1), we have
\[-(\rho_a-1/2)^2=\rho_a(1-\rho_a)-1/4\ge\rho_{a^\prime}(1-\rho_{a^\prime})
-1/4=-(\rho_{a^\prime}-1/2)^2.\]
Thus $\rho_a\le\rho_{a^\ast}$ and by Theorem \ref{th0.2}, 2), it follows that 
$\lim_{a\to\infty}\rho_a=1$. This proves the second part of 2).
 
Now consider the eigenvalue $1/4$. Let
\[\xi(s)=s(s-1)\Lambda(s).\]
Then  we get
\[c^\prime(1/2)=-\frac{c^\prime(1/2)}{c(1/2)}=-4\lim_{s\to 0}
\left(\frac{\Lambda^\prime(s)}{\Lambda(s)}+\frac{1}{s}\right)=-4\left(1+
\frac{\xi^\prime(0)}{\xi(0)}\right).\]
By \cite[pp.80-81]{Da} it follows that
\[2\left(1+\frac{\xi^\prime(0)}{\xi(0)}\right)=-\gamma+\log 4\pi,\]
where $\gamma$ denotes Euler's constant. Let $a^\ast:=4\pi e^{-\gamma}$. 
Then we get
\begin{equation}\label{2.3}
c^\prime(1/2)=-2\log a^\ast.
\end{equation}
By Theorem 0.2, 4), $1/4$ is an eigenvalue 
of $\Delta_a$ if and only if $a=a^\ast$.
Since all  eigenvalues $\mu_i(a)$ of $\Delta_a$ are non-increasing
functions of $a$ and have multiplicity one, it follows that
$\mu_0(a^\ast)=1/4$, $\mu_0(a)<1/4$, if $a>a^\ast$, and $\mu_0(a)>1/4$, if
$1\le a<a^\ast$. Since $\mu_0(a)=\rho_a(1-\rho_a)$, this proves 1) and the
first part of 2). 

It remains to show that all zeros $\rho\not=1/4$ are simple. Let $\rho$ be
a zeros of $\varphi_0(a,s)$ of multiplicity $\ge 2$. Then $\rho$ is also a 
zero of
$1+c(s)a^{1-2s}$ of multiplicity $\ge 2$. Thus
\[0=\frac{d}{ds}\left(1+c(s)a^{1-2s}\right)\big|_{s=\rho}=a^{1-2s}
\left(c^\prime(\rho)-2c(\rho)\log a \right).\]
and hence,  we get
\begin{equation}\label{2.4}
2\log a-\frac{c^\prime}{c}(\rho)=0.
\end{equation}
  
Let $\chi_{[a,\infty)}$ be the characteristic function of the interval
$[a,\infty)$ 
in $\R^+$. For $a\ge 1$ let
\[\hat E(z,s):=E(z,s)-\chi_{[a,\infty)}(y^s+c(s)y^{1-s})\]
be the truncated Eisenstein series. Then $\hat E(z,s)$ is square integrable
and its norm is given by the Maass-Selberg relation. 
For $r\in\R$, $r\not=0$, we have \cite[(7.42')]{Se}
\begin{equation}
\parallel \hat E(1/2+ir)\parallel^2=2\log a-\frac{c^\prime}{c}(1/2+ir)+
\frac{1}{r}\Im\left(c(1/2+ir)a^{-2ir}\right).
\end{equation} 
Let $\rho=1/2+ir$, $r\not=0$, be a zero of $1+c(s)a^{1-2s}$ of multiplicity
$\ge 2$. Then 
$$c(1/2+ir)=-a^{2ir}.$$
 Thus $\Im\left(c(1/2+ir)a^{-2ir}\right)=0$.
This together with (\ref{2.4}) implies that $\hat E(1/2+ir)=0$. On the other
hand, the zeroth Fourier coefficient is not identically zero. Hence we arrive
at a contradiction.

The case of real zeros is treated similarly. Let $\sigma\in(1/2,1)$ be a 
zero of $1+c(s)a^{1-2s}$ of multiplicity $\ge 2$. In this case the 
Maass-Selberg relation is
\begin{equation}
\parallel \hat E(\sigma)\parallel^2=\frac{a^{2\sigma-1}
-|c(\sigma)|^2a^{1-2\sigma}}{2\sigma-1}+2c(\sigma)\log a- c^\prime(\sigma)
\end{equation} 
\cite[(7.42'')]{Se}.
Since $\sigma$ is a zero, we have $|c(\sigma)|=a^{2\sigma-1}$. 
This together with (\ref{2.4}) implies
$\hat E(\sigma)=0$, which again leads to a contradiction. 

Finally note that by (\ref{2.2}) and (\ref{2.3}), $1/2$ is a double zero 
of $\varphi_0(a^\ast,s)$.  This completes the proof of Theorem 0.1.
\hfill$\Box$

We note that Th\'eor\`em 5 of \cite{CV} has an extension to hyperbolic surfaces
with  several cusps similar to Theorem \ref{th1.4}. In the case of $m$ cusps
the scattering matrix $C(s)$ is  an $m\times m$ matrix. In this case we study
the zeros of the function
\begin{equation}\label{2.5}
\phi(a,s):=\det(a^s+a^{1-s}C(s)).
\end{equation}
For a congruence subgroup
$\Gamma\subset \PSL(2,\Z)$ the entries of $C(s)$ consist of ratios of Dirichlet
$L$-series \cite{Hu}. 
If we multiply (\ref{2.5}) by the common denominator, we get a certain linear
 combination 
of products of Dirichlet $L$-functions.  As example, consider the case of the
 Hecke group 
$\Gamma=\Gamma_0(n)$, where $n=p_1\cdots p_r$ is a product of pairwise distinct
prime numbers $p_i$. Then the number of cusps of $\Gamma_0(n)\ba \Hb$ is 
$m=2^r$. The scattering matrix can be described in the following way. 
For a prime $p$ let
\[N_p(s):=\frac{1}{p^{2s}-1}\begin{pmatrix}p-1& p^s-p^{1-s}\\
p^s-p^{1-s}& p-1
\end{pmatrix}
\]
Then 
\[C(s)=\frac{\Lambda(2s-1)}{\Lambda(2s)}\left(N_{p_1}(s)\otimes\cdots 
\otimes N_{p_r}(s)\right)\]
\cite{He}. Suppose that  $n=p$ is a prime number. Then we get
\begin{equation}
\begin{split}
\det&\left(a^s\Id+a^{1-s}C(s)\right)\\
&=a^{2s}+2a\frac{p-1}{p^{2s}-1}\cdot
\frac{\Lambda(2s-1)}{\Lambda(2s)}+a^{2-2s}\frac{p^{2-2s}-1}{p^{2s}-1}\cdot
\frac{\Lambda(2s-1)^2}{\Lambda(2s)^2}\\
&=\frac{\prod_{k=1}^2\left(a^s(p^s+(-1)^k)\Lambda(2s)
+a^{1-s}(p^{1-s}+(-1)^k)\Lambda(2s-1)\right)}{(p^{2s}-1)\Lambda(2s)^2}. 
\end{split}
\end{equation}
Note that the numerator and the denominator have no common zeros. Therefore the
zeros of the function $\phi(a,s)$ coincide with the zeros
of the numerator. Let $a\ge 1$. Then it follows as above that the map
$s\mapsto s(1-s)$ induces a bijection between the zeros $\rho\not=1/2$ of 
$\phi(a,s)$ and the eigenvalues $\mu_j(a)\not=1/4$ of $\Delta_a$.
This implies that for $a\ge 1$ the complex zeros of
\[a^s(p^s\pm 1)\Lambda(2s)+a^{1-s}(p^{1-s}\pm 1)\Lambda(2s-1)\] 
lie all on the line $\Re(s)=1/2$. 
As for the real zeros, we observe that $\Gamma_0(p)$ has no residual
eigenvalues other than 0 \cite[Theorem 11.3]{Iw}. By the mini-max principle it
follows that 
$\mu_1(a)\ge 1/4$. Therefore, there exist at most two real zeros, $\rho_a$
and $1-\rho_a$, with $\rho_a\in(1/2,1)$. Moreover $\lim_{a\to\infty}\rho_a=1$.

\section{Imaginary quadratic fields}
\setcounter{equation}{0} 

In this section we derive Theorem \ref{th0.3} from  Theorem \ref{th1.5}. 

Let $K=\Q(\sqrt{-D})$ be an imaginary quadratic field, where $D$ is a square
free positive integer, 
and let $\cO_K$ be the ring of integers in $K$. Let 
$\Gamma_K=\PSL(2,\cO_K)$. 
Then $\Gamma_K\subset\PSL(2,\C)$ is a discrete subgroup of finite co-volume. The
number of cusps of $\Gamma_K\ba\Hb^3$ equals the class number $h_K$ of 
$K$. 
Assume that $h_K=1$. This is the case if and only if 
$D=1,2,3,7,11,19,43,67,163$ \cite{St}. Then $\Gamma_K\ba\Hb^3$ has a single 
cusps  which
we choose to be at $\kappa=\infty$. A fundamental domain $\cF_K$ of 
$\Gamma_K$ can be constructed as in \cite[Section 7.3]{EGM}, 
\cite[pp. 17-19]{Sw}. 
Let $F_K$ be a fundamental domain for the group
\[\mathrm{PSL}(2,\cO_K)_\infty=\left\{\begin{pmatrix}\alpha&\beta\\0&\lambda
\end{pmatrix}\in\mathrm{PSL}(2,\cO_K)\right\}.\]
If
$\mu,\lambda\in\cO_K$ 
generate $\cO_K$ and $\mu\not=0$, let 
\[H(\mu,\lambda):=\{(z,r)\in\Hb^3\colon |\mu z+\lambda|^2+|\mu|^2r^2\ge 1\}.\]
This is the set of points $w\in\Hb^3$ which lie above or on the hemisphere
$S(\mu,\lambda)$ in $\Hb^3$ which is  the set of all $(z,r)$ such that
$|\mu z+\lambda|^2+|\mu|^2 r^2=1$. This hemisphere has center $(\lambda/\mu,0)$
and radius $1/|\mu|$. Let 
\[B_K:=\bigcap H(\mu,\lambda),\]
where the intersection is taken over all $\mu,\lambda\in\cO_K$ as above. Then
a fundamental domain $\cF_K$ of $\Gamma_K$ is given by
\[\cF_K=\{(z,r)\in B_K\colon z\in F_K\}\]
\cite[Theorem 7.3.4]{EGM}, \cite[Section 3]{Sw}. Thus $\cF_K$ is the complement
in $F_K\times \R^+$ of Euclidean semi-balls with center at $(\lambda/\mu,0)$ and
radius $1/|\mu|$, $\mu\not=0$. Since $\mu\in\cO_K$, it follows that 
$|\mu|\ge 1$ and therefore, we have
\begin{equation}\label{3.2}
F_K\times [1,\infty)\subset \cF_K.
\end{equation}
For $\Q(i)$ and $\Q(\sqrt{-3})$,
the fundamental domain is a hyperbolic pyramid with one vertex at $\infty$.

Let $d_K$ be the discriminant of $K$ and
let $\zeta_K(s)$ be the Dedekind zeta function of $K$. Let
\[\Lambda_K(s)=\left(\frac{\sqrt{|d_K|}}{2\pi}\right)^s\Gamma(s)\zeta_K(s)\]
be the completed zeta function. 
Then the constant term of the Eisenstein
series $E(w,s)$ is given by
\begin{equation}\label{3.3}
E_0(y,s)=y^s+\frac{\Lambda_K(s-1)}{\Lambda_K(s)}y^{2-s}
\end{equation}
\cite[Section 8.3]{EGM}, \cite{ES}, and therefore, the scattering matrix
 $c_K(s)$ equals
\begin{equation}\label{3.4}
c_K(s)=\frac{\Lambda_K(s-1)}{\Lambda_K(s)}.
\end{equation} 
The completed zeta function satisfies
the functional equation $\Lambda_K(s)=\Lambda_K(1-s)$. All its zeros are
contained in the strip $0<\Re(s)<1$. This implies that $\Lambda_K(s-1)$
and $\Lambda_K(s)$ have no common zeros. Therefore for $a>0$, the zeros of 
\[\varphi_K(a,s)=a^s\Lambda_K(s-1)+a^{2-s}\Lambda_K(s)\]
coincide with the zeros of $E_0(a,s)$. 

We can now apply Theorem \ref{th1.5}. By (\ref{3.2}) we can take $b=1$.
Let $a\ge 1$ and let $\rho\not=1$ be a zero of $\psi_K(a,s)$. Then by
Theorem  \ref{th1.5}, 2), $\rho(2-\rho)$ is an eigenvalue of $\Delta_a$. This
implies 
that if $\rho(2-\rho)\ge 1$, then $\rho=1+ir$ with $r\in\R$, $r\not=0$. 
If $\rho(2-\rho)< 1$, then $\rho$ is real. Thus all complex zeros of 
$\varphi_K(a,s)$ are on the line $\Re(s)=1$.  

For the real zeros we need to consider the non-cuspidal eigenvalues 
 of $\Delta$. By the spectral resolution of the Laplacian, these eigenvalues 
are in 
one-to-one correspondence with the poles of the Eisenstein series $E(w,s)$ in 
the interval $(1,2]$ \cite[Proposition 6.2.2]{EGM}.
 On the other hand, using the
Maass-Selberg relations, it follows that
the poles of $E(w,s)$ in $(1,2]$ coincide with the poles of the scattering
matrix $c_K(s)$ in $(1,2]$. Since $\Lambda_K(s)$ has no zeros in $\Re(s)>1$
and the only pole in $\Re(s)>0$ is a simple pole at $s=1$, it follows from
(\ref{3.4}) that the only pole of $c_K(s)$ in $\Re(s)>1$ is a simple pole at 
$s=2$ which corresponds to the eigenvalue 0. This shows that $\Delta$ has no
non-cuspidal eigenvalues in the interval $(0,1)$.  
 By Theorem \ref{th1.5}, 4), it follows that  $\varphi_K(a,s)$ has  at most 
two real zeros $\rho_a$ and $2-\rho_a$ with
$\rho_a\in (1,2)$. Let $a<a^\prime$. Then by Theorem \ref{th1.5}, 1), we 
have
\[-(\rho_a-1)^2=\rho_a(2-\rho_a)-1\ge \rho_{a^\prime}(2-\rho_{a^\prime})-1=
-(\rho_{a^\prime}-1)^2.\]
This implies $\rho_a\le\rho_{a^\prime}$, and by 4) it follows that
$\lim_{a\to\infty}\rho_a=1$. 

Next we consider the eigenvalue $1$. Recall that $\Lambda_K(s)$ has simple
poles at $s=0$ and $s=1$. 
Using the functional equation, we get ${\mathrm Res}_{s=1}\Lambda_K(s)=-
{\mathrm Res}_{s=0}\Lambda_K(s)$. This implies $c_K(1)=-1$.
By Theorem \ref{th1.5}, 3), $1$ occurs as an eigenvalue $\mu_j(a)$ of
$\Delta_a$ if and only if $c_K^\prime(1)=-2\log a$. Let
\[\xi_K(s)=s(s-1)\Lambda_K(s).\]
Then
\[c_K^\prime(1)=-\frac{c_K^\prime(1)}{c_K(1)}=
\lim_{s\to 1}\left(\frac{\Lambda_K^\prime(s)}{\Lambda_K(s)}
-\frac{\Lambda_K^\prime(s-1)}{\Lambda_K(s-1)}\right)
=-2\left(1+\frac{\xi_K^\prime(0)}{\xi_K(0)}\right).\]
Let 
\[a^\ast_K:=\exp\left(1+\frac{\xi_K^\prime(0)}{\xi_K(0)}\right).\]
Then by the above, there is $j$ such that  $\mu_j(a)=1$ if and only if 
$a=a^\ast_K$. Assume that $a^\ast_K\ge 1$. Using that the eigenvalues $\mu_j(a)$
of $\Delta_a$ are non-increasing and have multiplicity one, we get 
$\mu_0(a^\ast_K)=1$, $\mu_0(a)<1$, if $a>a^\ast_K$, and $\mu_0(a)>1$, if
$a<a^\ast_K$. 
If $a^\ast_K<1$, then $\mu_0(a)<1$ for all $1\le a$.
Using that $\mu_0(a)=\rho_a(2-\rho_a)$, we get 1) and 2) of Theorem 
\ref{th0.3}. The simplicity of the zeros $\rho\not=1$ follows again from
the Maass-Selberg relations \cite[p.270]{EGM} in the same way as in the case 
of $\PSL(2,\Z)$.
\hfill$\Box$

The constant $a^\ast_K$  in Theorem \ref{th0.3} can be computed as follows. 
Let $w_K$ be the order of the group of units of $\cO_K$. Let $1,\omega$ be 
a basis for $\cO_K$ as a $\Z$ module. 
Since $K$ has class number one, we have
\[\zeta_K(s)=\frac{1}{w_K}{\sum}^\prime_{m,n}\frac{1}{|m+n\omega|^{2s}},\] 
the sum being taken for all integers $(m,n)\not=(0,0)$. 
Using the Kronecker limit formula \cite[p. 273]{La}, we get for 
$K=\Q(\sqrt{-D})$
\begin{equation}
a^*_{K}=\frac{e^{2+\gamma}}{2\pi\sqrt{D}}\begin{cases}
\frac{1}{2|\eta(i\sqrt{D})|^4},&\quad D\equiv 2,3 \mod 4\\ 
\frac{1}{\big|\eta\left(\frac{1+i\sqrt{D}}{2}\right)\big|^4},&\quad D\equiv 1
\mod 4, \end{cases}
\end{equation}
where $\gamma$ denotes Euler's constant and $\eta(z)$ is the Dedekind eta 
function. The values of the eta function  can be computed by the
Chowla-Selberg formula \cite[(2), p.110]{SC}.  Let $\Delta(z)=\eta(z)^{24}.$ 
Let $D\equiv 2,3 \mod 4$ and assume that $h(D)=1$. Then the Chowla-Selberg
formula gives 
\begin{equation} 
\Delta(\sqrt{D}i)=\frac{1}{(8\pi D)^6}\left\{\prod_{m=1}^{4D}
\Gamma\left(\frac{m}{4D}\right)^{\left(\frac{-4D}{m}\right)}\right\}^{3w_K}.
\end{equation}
So for $D=1$ we get
\[\Delta(i)=\frac{1}{(8\pi)^6}\left(\frac{\Gamma\left(\frac{1}{4}\right)}
{\Gamma\left(\frac{3}{4}\right)}\right)^{12}=
\frac{\Gamma\left(\frac{1}{4}\right)^{24}}{2^6(2\pi)^{18}},\] 
where in the right equality we used the duplication formula of the Gamma
function. This implies
$$a^\ast_{\Q(i)}=\frac{4\pi^2 e^{2+\gamma}}{\Gamma(1/4)^4}\approx
3.0068.$$
Similarly for $D=2$ we get 
\[\Delta(\sqrt{2}i)=\frac{1}{(2\pi)^68^6}
\left(\frac{\Gamma\left(\frac{1}{8}\right)\Gamma\left(\frac{3}{8}\right)}
{\Gamma\left(\frac{5}{8}\right)\Gamma\left(\frac{7}{8}\right)}\right)^6
=\frac{\left(\Gamma\left(\frac{1}{8}\right)\Gamma\left(\frac{3}{8}\right)\right)^{12}}{2^{15}(2\pi)^{18}},\]
and therefore
\[a^\ast_{\Q(\sqrt{-2})}=\frac{8\pi^2e^{2+\gamma}}
{\left(\Gamma\left(\frac{1}{8}\right)
\Gamma\left(\frac{3}{8}\right)\right)^{2}}\approx 
3.2581 .\]

\end{document}